\documentclass{amsart}

\usepackage{amsmath,amssymb,amsthm,mathtools}
\usepackage[hidelinks]{hyperref}
\newtheorem*{theorem}{Theorem}
\newcommand{\D}{\mathbb D}
\newcommand{\C}{\mathbb C}
\newcommand{\sgn}{\operatorname{sgn}}

\title{
Brennan's Conjecture and the Universal Integral Means Spectrum}

\author{Ilia Binder}
\address{Department of Mathematical and Computational Sciences,
University of Toronto,
Toronto, Ontario, Canada}
\email{ilia@math.utoronto.ca}

\author{Oleg Ivrii}
\address{Department of Theoretical Mathematics,
Tel Aviv University, Tel Aviv, Israel}
\email{ivrii@tauex.tau.ac.il}

\author{Istv\'an Prause}
\address{Mathematics, Faculty of Science and Engineering,
{\AA}bo Akademi University, Turku, Finland}
\email{istvan.prause@abo.fi}

\author{Stanislav Smirnov}
\address{Section de math\'ematiques,
Universit\'e de Gen\`eve, Geneva, Switzerland}
\email{stanislav.smirnov@unige.ch}

\begin{document}

\maketitle

\begin{abstract}
In October 2026, OpenAI released two preprints \cite{openai-brennan,openai-kraetzer},
establishing Brennan’s celebrated conjecture~\cite{brennan,brennan2} and disproving Kraetzer’s more general $t^2/4$ conjecture \cite{kraetzer, kraetzer-thesis}.
While validating these results,
we observed that much of their argument is redundant
and arrived at a drastically shorter proof of
a significantly stronger statement.
The goal of this note is to give a short and concise presentation of the reborn proof, sacrificing some of the context and consequences for brevity, and relegating them to a follow-up work.
In this particular case the human-AI interaction was fruitful and quite satisfactory, and hoping
for a similar future in other subjects, we would paraphrase Rick Blaine:
{\em ``AI, we hope this is the beginning of a beautiful friendship.''}
\end{abstract}

\section{Introduction}

\subsection{Brennan and we}

In 1978, Brennan~\cite{brennan,brennan2} conjectured that for every conformal
map $g:\Omega\to\mathbb D$ from a simply connected proper plane
domain, the area integral satisfies
\begin{equation}\label{eq:brennan}
\int_\Omega |g'(w)|^p\,dA(w)<\infty,
\qquad \frac43<p<4.
\end{equation}
Near an inward pointing slit ending at $0$, we would have $g(z)-g(0)\asymp \sqrt{z},~g'\asymp1/\sqrt{z}$,
so clearly the integral converges up to $p=4-$. Informally, the conjecture states that
we cannot create a more extremal example
by adding more slits as they start screening each other.

Together with many important questions in complex analysis,
Brennan's conjecture reduces to the study of
the integral means spectrum of conformal maps, or the multifractal spectra of harmonic measure,
as shown by Carleson and  Makarov  in \cite{carleson-makarov}.
For a univalent map of a unit disk $f\colon\D\to\C$, we define the beta spectrum by
\begin{equation}
M_t(f,r)=\int_0^{2\pi}|f'(re^{i\theta})|^t\,d\theta,
\qquad
\beta_f(t)=\limsup_{r\uparrow1}
\frac{\log M_t(f,r)}{\log(1/(1-r))}.
\label{eq:def-spectrum}
\end{equation}
We define the universal spectra by
\begin{equation}
\begin{aligned}
B(t)&=\sup\{\beta_f(t):~f\text{~univalent~in~}\D\},
\\
B_b(t)&=\sup\{\beta_f(t):~f\text{~univalent,~bounded~in~}\D\}.
\end{aligned}
\label{eq:def-univ}
\end{equation}
For instance, Brennan's conjecture reduces to the statement $B(-2)=1$,
which can be seen by taking $f=g^{-1}$, and observing
that the change of variables replaces $|g'|^p$ by $|f'|^{2-p}$.

The study of the beta spectra played a central role in geometric function theory
over the last decades, see the surveys~\cite{makarov,beliaev-smirnov-survey,hedenmalm-sola-survey}.
A general theory was developed, with the bounded spectra playing the central role,
as most other can be rewritten in its terms, for instance Makarov's identity \cite[Theorem 5.4]{makarov} states that
\begin{equation}
B(t)=\max\{B_b(t),3t-1\},
\label{eq:betab}
\end{equation}
A beautiful generalization of Brennan's conjecture was suggested by Kraetzer in his numerical work~\cite{kraetzer,kraetzer-thesis}:
\begin{equation}
B_b(t)=t^2/4,~t\in[-2,2],\quad B_b(t)=|t|-1,~t\notin[-2,2],
\label{eq:t24}
\end{equation}
and the authors and their colleagues invested significant time in laying the foundations for it.
Although many results supported the Kraetzer conjecture
\cite{carleson-jones,carleson-makarov,jones-makarov,beliaev-smirnov,beliaev},
recently evidence against it began to accumulate
\cite{ivrii-quasicircles,hedenmalm-zero-packing,
astala-ivrii-perala-prause,ivrii-makarov}.

\subsection{Brennan and AI}

Two  OpenAI preprints released in October 2026 presented a dramatic advance.
The Brennan conjecture was proved, while Kraetzer's conjecture
\eqref{eq:t24} was shown to be false at $t=-1$.
The proofs in \cite{openai-brennan,openai-kraetzer} rely on a family of potentials whose logarithms are superharmonic. A critical-point map records the target associated with each critical point, and the sign of its Jacobian distinguishes maxima from saddles. Comparing their heights and integrating over targets gives a weighted signed-area inequality.
Most ingredients are familiar, and have appeared before. The critical maps have antecedents in quasiconformal reflection and critical points of conformal metrics \cite{ahlfors-reflections,chuaqui-osgood94}. Morse theory and, more recently, zero-dimensional persistence naturally describe the topological bookkeeping \cite{milnor-morse,persistence}. Area inequalities for integral-means estimates also have a long history, particularly in Hedenmalm and Shimorin's weighted Bergman-space approach \cite{hedenmalm-shimorin}.

To prove a strict inequality for $B(-1)$, OpenAI's argument adds a probabilistic construction. A positive growth operator encodes derivative moments at successive scales, and a derivative-weighted affine-covariant probability law allows one to average the geometric comparison. The strict inequality is guaranteed by a global harmonic obstruction. The operator approach has predecessors in thermodynamic formalism \cite{ruelle-transfer}; a closer conformal-mapping analog is the positive-operator method for random conformal snowflakes \cite{beliaev,beliaev-smirnov}. The construction of the law parallels invariant-measure averaging \cite{krylov-bogolyubov} and averaged transport with transport of stationary random measures \cite{last-thorisson,gentner-last}.

We want to emphasize that while none of the techniques appear unprecedented, their combination in this context was  missed by the human experts.
On the other hand, the AI proofs contain a significant number of derivations that, while correct,
are redundant and do not seem to contribute either to strength or clarity of the results.

\subsection{AI and we}

When the OpenAI proofs appeared, at first we were puzzled
by a seemingly mysterious chain of thought.
However, within a day we independently
realized that many of the arguments were redundant,
and with a shorter proof we arrive at a stronger result.
Only some of us used AI agents to help with understanding,
and the write up below was done without AI help.

It is still unclear whether the redundant parts we discarded can be of some use elsewhere. Even the remaining argument, while short, is somewhat mysterious.
It clearly can give better estimates, but it is not clear whether
it can ever give a sharp estimate from above
or provide description of the extremal map.

At present, it seems that further advances in this direction would be computationally
intensive, and we have tried to use AI tools to this effect.
The results will be reported in a separate preprint.

Overall, this was a pleasant and fruitful AI-human interaction in mathematics:
there was a spectacular progress,
human community was brought together,
new results are about to follow, and
no robots or humans were harmed.

\subsection{Results}

\begin{theorem}
The bounded universal spectrum
satisfies
$B_b(t)\le t^2/4$ for $t\in[-2,2]$,
and the inequality is strict for $t\in]-2,2[\setminus\{0\}$.
\end{theorem}

From Makarov's identity \cite[Theorem 5.4]{makarov},
$B(t)=\max\{B_b(t),3t-1\},$
so it immediately follows
that
$B(t)< t^2/4$ for $t\in]-2,6-4\sqrt2[\setminus\{0\}$, since $3t-1\le t^2/4$
throughout this range.

For clarity, below we give a full argument only for a non-strict inequality.
One can obtain a strict estimate on $]-2,2[\setminus\{0\}$
by perturbing the weight below to e.g.
\begin{equation}
V^*_\xi(z)=\frac{(1-|z|^2)^{k^*}}
{|z|^{k^*}(f(z)-\xi)}
~e^{
\varepsilon^*\left(
\arctan \left|(1-|z|^2)z\,\frac{f''(z)}{f'(z)}-\frac t2\right|^2
-\frac\pi4
\right)},
\label{eq:newv}
\end{equation}
and replacing the pointwise positivity by an averaged one.
We defer the details to a separate paper.

We begin with the introduction of a potential $V_\xi$ and auxiliary function $G$ \eqref{eq:VGT}. Jacobian of $G$ distinguishes maxima and saddles of the potential. Subharmonicity of $\log V_\xi$ and the index theorem imply that there are  more maxima than saddles above every level. By integrating this identity with respect to height,  we obtain the signed area inequality \eqref{eq:signed-area}. By integration by parts, we obtain an expression for integral means in terms of signed area with a controlled boundary term correction \eqref{eq:radial}. There is a singularity near the origin for $t>0$. It is dealt with by localizing the weight in \eqref{eq:signed-area} near the unit circle. All of these estimates are done under the assumption that our conformal map can be extended to be analytic in a neighborhood of the closed unit disk. This is not restrictive, since all our estimates on integral means are uniform for a fixed radius $r$.

\section{Proofs}
For a map $\phi$ bounded and univalent in $\D$ and
$1/2<\rho<1$, set $f(z)=\phi(\rho z)/\rho$.
Then $f$ is univalent and bounded in a neighborhood of $\overline\D$, and we let
$\rho\uparrow1$ at the end.
Fix $t\ne0$ and set
\begin{equation}
q=2-t>0,\qquad k>\frac q4,\qquad s(z)=1-|z|^2,
\qquad b=(k-1)q>-1,
\label{eq:parameters}
\end{equation}
the latter inequality holds because
$(k-1)q+1>(\frac q4-1)q+1=(\frac q2-1)^2\ge0$.

\subsection{Auxiliary functions and their basic properties}
For $\xi\in\C\setminus f(\partial\D)$, define
\begin{equation}
V_\xi(z)=\frac{s(z)^k}{|f(z)-\xi|},\qquad
G(z)=f(z)+\frac{s(z)f'(z)}{2k\bar z},\qquad
T(z)=\frac{|z|s(z)^{k-1}}{|f'(z)|}.
\label{eq:VGT}
\end{equation}
Writing $\ell_\xi=\log V_\xi$, away from the possible pole $z_\xi: ~ f(z_\xi)=\xi$ we have strict superharmonicity:
\begin{equation}
\Delta\ell_\xi=\Delta (k \log s-\log |f(z)-\xi|) = k\Delta \log s=-{4k}{s^{-2}}<0.
\label{eq:superharmonic}
\end{equation}
Also $V_\xi\to0$ at
$\partial\D$, hence every positive superlevel set is compactly contained
in $\D$.
Similarly
\begin{equation}
\partial_z\ell_\xi
=\partial_z (k \log s-\log |f(z)-\xi|)
=-\frac{k\bar z}{s}
-\frac12\frac{f'}{f-\xi}
=-(G-\xi)\frac{k\bar z}{s(f-\xi)}
~,
\label{eq:dell}
\end{equation}
so a point $z\ne0,z_\xi$ is critical for $V_\xi$ or $\ell_\xi$ exactly when
\begin{equation}
G(z)=\xi.
\label{eq:critical-map}
\end{equation}
Note that $G$ is defined on $\D\setminus\{0\}$. If $\xi\neq f(0)$, then there is no critical point omitted at the origin.

At a critical point of $V_\xi$, we have
\begin{equation}
|f(z)-\xi|=\left|G-\xi-\frac{s f'}{2k\bar z}\right|=\frac{s|f'|}{2k|z|},~
\quad V_\xi(z)=
\frac{s^k}{|f-\xi|}
=
\frac{s^k2k|z|}{|sf'|}
=2kT(z).
\label{eq:critical-height}
\end{equation}
Differentiating $G$ we calculate
\begin{eqnarray}
\partial_z G=
f'+\frac1{2k\bar z}
\partial_z\left(sf'\right)
=\frac{f'}{2k}\left(2k-1+\frac{sf''}{\bar z f'}
\right)
,~\\
\partial_{\bar z} G=
\frac{f'}{2k}\partial_{\bar z}\left(\frac{s}{\bar z}\right)=
-\frac{f'}{2k\bar z^2}
,
\label{eq:difg}
\end{eqnarray}
hence we can expand
\begin{equation}
\det DG=
|\partial_z G|^2-|\partial_{\bar z} G|^2
=\frac{|f'|^2}{4k^2}
\left(\left|2k-1+\frac{s f''}{\bar z f'}\right|^2-\frac1{|z|^4}\right).
\label{eq:jacobian}
\end{equation}
Differentiating \eqref{eq:dell} -- a somewhat tedious calculation -- we get a simplified result
at a critical point where $G(z)=\xi$
by using \eqref{eq:critical-height}:
\begin{equation}
\det D^2\ell_\xi=
-\det D(2\partial_z\ell_\xi)=-\frac{4k^2|z|^2}{s^2|f-\xi|^2}\det DG=-\frac{16k^4|z|^4}{s^4|f'|^2}\,\det DG.
\label{eq:hessian-jacobian}
\end{equation}
By \eqref{eq:superharmonic} $\ell_\xi$ is superharmonic, so $\operatorname{tr}D^2\ell_\xi<0$; hence
$\det DG>0$ corresponds to a saddle and $\det DG<0$ to a local maximum.
Note that $$
G-f\longrightarrow0\quad\text{at }\partial\D,
\qquad
|G(z)|\longrightarrow\infty\quad\text{as }z\to0.
$$ For almost every $\xi$, Sard's theorem and compactness give finitely many,
all non-degenerate, solutions of \eqref{eq:critical-map}.

\subsection{Counting critical points}

Suppose that the height $h$ is regular, meaning that the corresponding level set of $V_\xi$ contains no critical points. Then every component
$U_j$ of $\{V_\xi>h\}$ is a Jordan domain: a hole would force an
interior minimum of $\ell_\xi$, contradicting the superharmonicity in
\eqref{eq:superharmonic}.
On $\partial U_j$, the vector field $\nabla\ell_\xi$ points inward.
The index theorem gives
\begin{equation}
1=M_j-S_j+P_j,
\label{eq:index}
\end{equation}
where $M_j,S_j,P_j$ count maxima, saddles, and poles in $U_j$.
Since $f$ is univalent, there is at most one pole altogether; summing
\eqref{eq:index} over the components gives
\begin{equation}
M(h)-S(h)=\#\{U_j\}-P(h)\ge0,
\label{eq:count}
\end{equation}
where $M(h),S(h),P(h)$ count maxima, saddles, and poles above the regular height $h>0$.

Integrating \eqref{eq:count} over heights and using \eqref{eq:critical-height}, we obtain
$$\sum_{G(z)=\xi}T(z)^q\sgn\det DG(z)
=-\frac{q}{(2k)^q}\int_0^\infty h^{q-1}\bigl(M(h)-S(h)\bigr)\,dh\le0.$$
Apply the real area formula separately on
$\{\det DG>0\}$ and $\{\det DG<0\}$.
For every Borel function $0\le\psi\le1$ satisfying
$$
\int_{\D}\psi(G(z))T(z)^q|\det DG(z)|\,dA(z)<\infty,
$$
subtracting the two formulas
gives
\begin{multline}
\int_{\D}\psi(G(z))T(z)^q\det DG(z)\,dA(z)\\
=\int_{\mathbb C}\psi(\xi)
\sum_{G(z)=\xi}
T(z)^q\,\operatorname{sgn}\det DG(z)\,dA(\xi)
\le0.
\label{eq:signed-area}
\end{multline}

\subsection{Conclusion}

Denote
\begin{equation}
M(r)=M_t(f,r)=\int_0^{2\pi}|f'(re^{i\theta})|^t\,d\theta,\qquad H(re^{i\theta})=|f'(re^{i\theta})|^t.
\label{eq:MH}
\end{equation}
Since $f'\neq0$, differentiating in $r$ and
using
$\Delta |f'|^t=t^2|f'|^t\left|\frac{f''}{f'}\right|^2$
we get
\begin{equation}
\begin{aligned}
M'(r)&=t\int_0^{2\pi}H\,\Re(e^{i\theta} f''/f')\,d\theta,
\\
M''(r)+\frac{M'(r)}r&=t^2\int_0^{2\pi}H|f''/f'|^2\,d\theta.
\end{aligned}
\label{eq:Mderivatives}
\end{equation}

Let
\begin{equation}
X(r)=\frac{r^{q-1}(1-r^2)^{b+2}}{t^2}M'(r)
+\left[
\frac{2k(t+2)}{t^2}r^q(1-r^2)^{b+1}
+\frac1t r^{q-2}(1-r^2)^{b+2}
\right]M(r).
\label{eq:X}
\end{equation}
Equations \eqref{eq:jacobian} and \eqref{eq:Mderivatives} give
$$\begin{aligned}4k^2\int_0^{2\pi}T^q\det DG\,r\,d\theta&=\frac{r^{q-1}s^{b+2}}{t^2}\left(M''+\frac{M'}r\right)\\
&\quad+\frac{2(2k-1)}t r^qs^{b+1}M'+\left((2k-1)^2r^{q+1}-r^{q-3}\right)s^bM\\
&=X'(r)+\frac{4k^2}{t^2}r^{q-1}s^b\left(4r^2-\frac qk\right)M.\end{aligned}$$
The last equality follows by differentiating \eqref{eq:X}, using $2-q=t$ and $(2k-1)t+b+2=k(t+2)$.
Since $M$ and $M'$ are bounded near $1$ and $b>-1$, we have $X(1^-)=0$; integrating from $r_0$ to $1$ gives
\begin{equation}
\int_{r_0<|z|<1}T^q\det DG\,dA
=-\frac{X(r_0)}{4k^2}
+\frac1{t^2}\int_{r_0}^1
r^{q-1}(1-r^2)^bM(r)
\left(4r^2-\frac qk\right)dr,
\label{eq:radial}
\end{equation}

Importantly, all interior terms in \eqref{eq:radial} involving $f''/f'$ have
been absorbed into the boundary term $X(r_0)$!

Choose $R_0\in(1/2,1)$ such that $4R_0^2>q/k$.
Then $4r^2-q/k$ is bounded below by a positive constant on $[R_0,1)$,
independently of $\rho$.
We claim that, uniformly in \(\rho\),
\begin{equation}
\int_{R_0}^1 r^{q-1}(1-r^2)^bM_t(f,r)\,dr\le C.
\label{eq:tail}
\end{equation}
The weight $\phi$ is bounded, hence Cauchy
estimate gives a disk $\Delta\subset\C$ and
$\varepsilon_0>0$, independent of $\rho$, such that
\begin{equation}
G(z)\in\Delta\quad(R_0\le|z|<1),
\qquad
G(z)\notin\Delta\quad(0<|z|<\varepsilon_0).
\label{eq:localization}
\end{equation}
Indeed $G$ is uniformly bounded for $|z|\ge R_0$, whereas
$G(z)=f'(0)/(2k\bar z)+O(1)$ as $z\to0$.

Since $b>-1$, the signed-area integral is absolutely convergent near the unit circle. We can now apply \eqref{eq:signed-area} with $\psi=\mathbf1_\Delta$.
Since the outer annulus maps into $\Delta$ and $G^{-1}(\Delta)$ avoids $|z|<\varepsilon_0$,
$$\int_{R_0<|z|<1}T^q\det DG\,dA\le\int_{\varepsilon_0\le|z|\le R_0}T^q|\det DG|\,dA\le C.$$
The last bound is uniform in $\rho$ by compactness.
Together with \eqref{eq:radial} and $X(R_0)=O(1)$, this proves \eqref{eq:tail}.

Now $f'(z)=\phi'(\rho z)$, so
\begin{equation}
M_t(f,r)=M_t(\phi,\rho r).
\label{eq:dilation}
\end{equation}
Moreover $M_t(\phi,r)$ is non-decreasing in $r$ since $|(\phi')|^t$ is subharmonic.
Letting $\rho\uparrow1$ in \eqref{eq:tail} and using monotone
convergence gives
\begin{equation}
\int_{R_0}^1(1-r)^bM_t(\phi,r)\,dr<\infty.
\label{eq:weighted}
\end{equation}
Hence, again by monotonicity,
\begin{equation}
M_t(\phi,r)\lesssim(1-r)^{-b-1},
\qquad
\beta_\phi(t)\le b+1.
\label{eq:beta}
\end{equation}
Finally $k\downarrow q/4$ gives
\begin{equation}
b+1=(k-1)q+1\downarrow
\left(\frac q4-1\right)q+1=\frac{t^2}{4}.
\label{eq:limit}
\end{equation}
So $\beta_\phi(t)\le t^2/4$ for every $-2\leq t<2$, $t\neq0$ and
every bounded univalent $\phi$.  The case $t=0$ is immediate.
For $t=2$, boundedness gives
$\int_0^1M_2(\phi,r)r\,dr=\operatorname{Area}(\phi(\D))<\infty$,
and by monotonicity $M_2(\phi,r)=O((1-r)^{-1})$, so
$\beta_\phi(2)\le1$.

\bibliographystyle{plain}
\bibliography{kraetzer_project}
\end{document}